\documentclass[letterpaper,10pt]{amsart}
\usepackage{amsmath,amssymb,amsxtra,amsthm}
\usepackage{graphicx}

\theoremstyle{plain}
\newtheorem{theorem}{Theorem}

\newtheorem{lemma}{Lemma}
\newtheorem{proposition}{Proposition}

\theoremstyle{definition}

\newtheorem{remark}{Remark}

\newcommand{\B}{\mathbb}
\newcommand{\C}{\mathcal}
\newcommand{\F}{\mathfrak}

\newcommand{\ga}{\alpha}
\newcommand{\gb}{\beta}
\newcommand{\gd}{\delta}

\newcommand{\gf}{\varphi}

\newcommand{\gl}{\lambda}

\newcommand{\gs}{\sigma}

\newcommand{\gD}{\Delta}
\newcommand{\gF}{\Phi}

\DeclareMathOperator{\md}{mod\ }

\newcommand{\leg}[2]{\left(\frac{#1}{#2}\right)}

\begin{document}
\title{Transcendence of the Gaussian Liouville number and relatives}
\author{Peter Borwein}
\author{Michael Coons}
\thanks{Research supported in part by grants from
 NSERC of Canada and MITACS.}%
\address{Department of Mathematics, Simon Fraser University, B.C., Canada V5A 1S6}
\email{pborwein@cecm.sfu.ca, mcoons@sfu.ca}
\subjclass[1991]{Primary 11J81; 11A05 Secondary 39B32; 11N64}%
\keywords{Transcendence of Power Series, Completely Multiplicative Functions}%
\date{\today}

\begin{abstract} {\em The Liouville number}, denoted $l$, is defined by $$l:=0.100101011101101111100\ldots,$$ where the $n$th bit is given by $\frac{1}{2}(1+\gl(n))$; here $\gl$ is the Liouville function for the parity of prime divisors of $n$. Presumably the Liouville number is transcendental, though at present, a proof is unattainable. Similarly, define {\em the Gaussian Liouville number} by $$\gamma:=0.110110011100100111011\ldots$$ where the $n$th bit reflects the parity of the number of rational Gaussian primes dividing $n$, $1$ for even and $0$ for odd. 

In this paper, we prove that the Gaussian Liouville number and its relatives are transcendental. One such relative is the number $$\sum_{k=0}^\infty\frac{2^{3^k}}{2^{3^k2}+2^{3^k}+1}=0.101100101101100100101\ldots,$$ where the $n$th bit is determined by the parity of the number of prime divisors that are equivalent to $2$ modulo $3$.

We use methods similar to that of Dekking's proof of the transcendence of the Thue--Morse number \cite{Dek1} as well as a theorem of Mahler's \cite{Mahl1}. (For completeness we provide proofs of all needed results.) This method involves proving the transcendence of formal power series arising as generating functions of completely multiplicative functions.



\end{abstract}

\maketitle


The {\em Liouville function} is the unique completely multiplicative function $\gl$ with the property that for each prime $p$, $\gl(p)=-1$. Denote the sequence of $\gl$ values by $\F{L}$.

Recall that a binary sequence is {\em simply normal} if each bit occurs with asymptotic frequency $\frac{1}{2}$, and {\em normal to base $2$} if each possible block of length $k$ occurs with asymptotic frequency $2^{-k}.$

The prime number theorem is equivalent to the simple normality of $\F{L}$; it is believed that $\F{L}$ is normal to base $2$, though a proof of normality is at present unattainable. 

We can get at some properties of $\F{L}$, though the sequence definitions of these are somewhat cumbersome. Define {\em the Liouville number} $l$ as $$l:=\sum_{n\in\B{N}}\left(\frac{1+\gl(n)}{2}\right)\frac{1}{2^n}=0.100101001100011100001\ldots.$$ The Liouville number is merely the translation of $\F{L}$ into ones and zeroes, and then a concatenation of the resulting sequence into a binary number. Properties of the number $l$ are properties of the sequence $\F{L}$. One noteworthy property is that $l$ is irrational.

\begin{theorem}\label{irrationall} If $l$ is the Liouville number, then $l\notin \B{Q}$.
\end{theorem}

\begin{proof} See proof of the general statement, Theorem \ref{irrationalgen}.
\end{proof}

The irrationality of $l$ tells us that the sequence $\F{L}$ is not eventually periodic. A fundamental question arises, which at present we are unable to answer. One assumes that $l$ is transcendental, though this is currently unapproachable. By considering other completely multiplicative functions like $\gl$, we may approach similar questions with success. 

We will consider the class of functions $f:\B{N}\to\{-1,1\}$, denoted $\C{F}(\{-1,1\})$ (this notation was used by Granville and Soundararajan in \cite{Gran1}), and call such $f$ {\em generalized Liouville functions}; this is not a particularly standard notation. Consider the following example.

Let $g$ be the completely multiplicative function defined on primes $p$ by $$g_p=\begin{cases} -1 & \mbox{if } p\equiv 3 (\md 4)\\ 1 & \mbox{otherwise}.\end{cases}$$ As the function $g$ takes the value $-1$ on those primes which are rational Gaussian primes, we call $g$ {\em the Gaussian Liouville function}. Denote by $\F{G}$ the sequence of values of $g$ and define $$\gamma:=\sum_{n\in\B{N}}\left(\frac{1+g_n}{2}\right)\frac{1}{2^n}$$ as the {\em Gaussian Liouville number}. One can show that $g_n=(-1/n)$ where $(\cdot/n)$ is the Jacobi symbol modulo $n$. The Gaussian Liouville number is easily seen to be irrational. Indeed, it is a corollary to the following generalization of Theorem \ref{irrationall}.

\begin{theorem}\label{irrationalgen} Suppose that $f:\B{N}\to\{-1,1\}$ is a completely multiplicative function, with $f(p)=-1$ for at least one prime $p$, and $\F{F}$ it's sequence of values. If $\gf:=\sum_{n\in\B{N}}\left(\frac{1+f(n)}{2}\right)\frac{1}{2^n}$, then $\gf\notin \B{Q}$.
\end{theorem}

\begin{proof} Towards a contradiction, suppose that $\gf\in \B{Q}$. Thus $\F{F}$ is eventually periodic, say the sequence is periodic after the $M$--th term and has period $k$. Now there is an $N\in\B{N}$ such that for all $n\geq N$, we have $nk>M$. Let $p$ be a prime for which $f(p)=-1$. Then $$f(pnk)=f(p)f(nk)= -f(nk).$$ But $pnk\equiv nk (\md k)$, a contradiction to the eventual $k$--periodicity of $\F{F}$. Hence $\gf\notin\B{Q}$.
\end{proof}

Though the transcendence of the Liouville number is unattainable, it is possible to establish the transcendence of the Gaussian Liouville number and many of its reatives. The proof of this result is contained in Section \ref{proofg}, and rests on the fact that the generating function of the sequence $\F{G}$ satisfies a useful functional equation (see Lemma \ref{func}). This functional equation leads to a striking power series representation of the functional equation, which is of interest (an example of such a series representation is given in Section \ref{2mod3}), and may lead to a quick transcendence result. 

As an example of the usefulness of such a representation, we prove in Section \ref{proofg} that the generating function $G(z)$ of the Gaussian Liouville sequence is $$G(z)=\sum_{k=0}^{\infty}\frac{z^{2^k}}{1+z^{2^{k+1}}}.$$ Noting that $\gamma=G(1/2)$, a quick transcendence result may be gained after some rearrangement and application the following theorem of Duverney \cite{Duv1}.

\begin{theorem}[Duverney, 2001] Let $a\in\B{N}\backslash\{0,1\}$ and let $b_n$ be a sequence of rational integers satisfying $|b_n|=O(\eta^{-2^n})$ for every $\eta\in(0,1).$ Suppose that $a^{2^n}+b_n\neq 0$ for every $n\in\B{N}.$ Then the number $$S=\sum_{n=0}^\infty\frac{1}{a^{2^n}+b_n}$$ is transcendental.
\end{theorem}

Though Duverney's theorem is very easy to apply\footnote{In Duverney's theorem, we need only set $a=2$, $b_n=2^{-2^n}$, and let $\eta\in(0,1)$. Then for all $n\in\B{N}$, $\eta<2$ implies that $|b_n|=2^{-2^n}<\eta^{-2^n}=O(\eta^{-2^n}).$}, it does not readily generalize enough for our purposes. With the goal of developing more general methods, we prove the transcendence of the Gaussian Liouville number without the application of this theorem. 

An ingredient is the proof of the transcendence of the generating function of the sequence $\F{G}$. While this may be proven from results found with difficultly in the literature (see Fatou \cite{Fat1}), for completeness and ease of reading, we offer a direct and elementary proof of this. In fact, one can prove transcendence of various functions using the following beautiful result.

\begin{theorem}[Fatou, 1906]\label{algtrans} A power series whose coefficients take only finitely many values is either rational or transcendental.
\end{theorem}

There is a recent proof of this by Allouche \cite{All1}; see also Borwein and Coons \cite{BC2} and Borwein, Erd\'elyi and Littmann \cite{BEL}.

The method used in our proof can be used to prove more general results regarding other completely multiplicative functions. Section \ref{gensection} contains these results (for an account of the properties of these functions see \cite{BCC1}).

A few historical remarks are in order. The irrationality of the values of power series similar to those of our investigation have been studied by, among others, Erd\H{o}s, Golomb, Mahler, and Schwarz. Erd\H{o}s \cite{Erd1} proved that the series $$\sum_{n=1}^\infty\frac{z^n}{1-z^n}=\sum_{n=1}^\infty d(n)z^n,$$ where $d(n)$ is the divisor counting function, gives irrational values at $z=\frac{1}{t}$ for $t=2,3,4,\ldots$, and Allouche \cite{All1} has shown this function to be transcendental, but all values are still open, for example $z=\frac{1}{2}$ presumably gives a transcendental value. Indeed Erd\H{o}s writes \cite{Erd2}, \begin{quote}``It is very annoying that I cannot prove that $\sum_{n=1}^\infty\frac{1}{2^n-3}$ and $\sum_{n=2}^\infty\frac{1}{n!-1}$ are both irrational (one of course expects that $\sum_{n=1}^\infty\frac{1}{2^n+t}$ and $\sum_{n=2}^\infty\frac{1}{n!+t}$ are irrational and in fact transcendental for every integer $t$.)''\end{quote} Partially answering Erd\H{o}s' question, Borwein \cite{Bor3} has shown that  for $q\in\B{Z}$ with $|q|>1$ and $c\in\B{Q}$ with $c\neq 0$ and $c\neq -q^n$ $(n\in\B{N})$, $$\sum_{n=1}^\infty \frac{1}{q^n+c}\quad\mbox{and}\quad\sum_{n=1}^\infty \frac{(-1)^n}{q^n+c}$$ are irrational; the special values $c=-1$ and $q=2$ give that the sum of the reciprocals of the Mersenne numbers is irrational. Later, Golomb proved in \cite{Gol1} that the values of the functions $$\sum_{n=0}^\infty\frac{z^{2^n}}{1+z^{2^n}}\qquad\mbox{and}\qquad\sum_{n=0}^\infty\frac{z^{2^n}}{1-z^{2^n}}$$ are irrational at $z=\frac{1}{t}$ for $t=2,3,4,\ldots,$ the interesting special case of which is that the sum of the inverse of the Fermat numbers is irrational; transcendence of the sum of the inverse of the Fermat numbers is implied by Duverney's theorem. Schwarz \cite{Sch1} has given results on series of the form $$\sum_{n=0}^\infty \frac{z^{k^n}}{1-z^{k^n}}.$$ In particular, he showed that this function is transcendental at certain rational values of $z$ when $k\geq 2$ is an integer. We take these results further and prove that for a completely multiplicative function $f$, with recursive relations $f_p=\pm1$ and $f_{pk+i}=f_i$ for $i=1,2,\ldots,p-1,$ the series $$\sum_{k=0}^\infty f(n)z^n=\sum_{k=0}^\infty \frac{f_p^k\Phi(z)}{1-z^{p^{k+1}}},$$ with $\Phi(z)=\sum_{i=1}^{p-1} f_iz^{p^ki}$, is transcendental (see Theorem \ref{gentransfunc} and Proposition \ref{Fseries}); it is interesting to note that when $f_i=(i/p)$ for $p\nmid i$ is the Legendre symbol, the polynomial $\Phi(z)$ is the $p$th degree Fekete polynomial. Patterns in the sequence of values of such $f$ have been studied by Hudson \cite{Hud1,Hud2}.

Mahler's results are too numerous to mention, and it seems likely that at least some of the historical results mentioned here were known to him as early as the 1920s (see \cite{Mahl2}). Mahler was one of the first to consider the links between functional equations and transcendence. The book \cite{Nish1} serves as a testimonial to this.

The methods used in this investigation were inspired by Michel Dekking's proof of the transcendence of the Thue--Morse number \cite{Dek1}, which is based on a method of Cobham \cite{Cob1}. 

\section{The Liouville function for primes 2 modulo 3}\label{2mod3}

As a striking example of a power series representation of a generating function consider the completely multiplicative function $t_n$ where $$t_3=1\quad\mbox{and}\quad t_p=\begin{cases}-1 & \mbox{if } p\equiv 2 (\md 3)\\ 1 & \mbox{if } p\equiv 1 (\md 3).\end{cases}$$ We have the relations $$t_{3n}=t_n,\quad t_{3n+1}=1,\quad\mbox{and}\quad t_{3n+2}=-1.$$ Denote the generating function of $(t_n)_{n\in\B{N}}$ as $T(z)=\sum_{n=1}^\infty t_nz^n$. Then $$T(z)=\sum_{n=1}^\infty t_{3n}z^{3n}+\sum_{n=0}^\infty t_{3n+1}z^{3n+1}+\sum_{n=0}^\infty t_{3n+2}z^{3n+2}
=T(z^3)+(z-z^2)\frac{1}{1-z^3},$$ which gives the following result.

\begin{lemma} If $T(z)=\sum_n t_nz^n$, then $$T(z^3)=T(z)-\frac{z}{1+z+z^2}.$$
\end{lemma}

Using this lemma we have $$T(z^{3^m})=T(z)-\sum_{k=0}^{m-1}\frac{z^{3^k}}{1+z^{3^k}+z^{3^k2}}.$$ Denote the sum by $U_m(z)$; that is, $$U_m(z)=\sum_{k=0}^{m-1}\frac{z^{3^k}}{1+z^{3^k}+z^{3^k2}}.$$ If $|z|<1$, taking the limit as $m\to\infty$, gives the desired series expression. 

\begin{proposition} If $|z|<1$, then the generating function of $(t_n)_{n\in\B{N}}$ has the closed form  $$T(z)=\sum_{k=0}^{\infty}\frac{z^{3^k}}{1+z^{3^k}+z^{3^k2}}.$$
\end{proposition} 

Application of the general results proved in Section \ref{gensection} gives the following result.

\begin{theorem} The function $T(z)$ is transcendental; furthermore, $T(\ga)$ is transcendental for all nonzero algebraic numbers $\ga$ with $|\ga|<1$.
\end{theorem}

\section{The Gaussian Liouville Function}\label{proofg}

As before, {\em the Gaussian Liouville function} $g$ is the completely multiplicative function defined on the primes by $$g_p=\begin{cases} -1 & \mbox{if } p\equiv 3 (\md 4)\\ 1 & \mbox{otherwise}.\end{cases}$$ Also,  denote by $\F{G}$ the sequence of values of $g$, and by $$\gamma:=\sum_{n\in\B{N}}\left(\frac{1+g_n}{2}\right)\frac{1}{2^n}$$ the {\em Gaussian Liouville number}.

The first few values of $g$ are $$\F{G}=(1,1,-1,1,1,-1,-1,1,1,1,-1,-1,1,-1,-1,1,1,1,-1,1,1,-1,\ldots).$$

Elementary observations tell us that the occurrence of primes that are 3 modulo 4 in prime factorizations are fairly predictable. One has the following implications: \begin{align*} n\equiv 1 (\md 4)\ &\Rightarrow\ g_n=1\\
n\equiv 3 (\md 4)\ &\Rightarrow\ g_n=-1\\
n\equiv 0 (\md 2)\ &\Rightarrow\ g_{2n}=g_n,\end{align*}
which give the recurrence relations for the sequence $\F{G}$ as
\begin{equation}\label{recs}g_1=1,\quad g_{2n}=g_n,\quad g_{4k+1}=-g_{4k+3}.\end{equation} This is not so surprising recalling that $g_n=(-1/n)$ where $(\cdot/n)$ is the Jacobi symbol modulo $n$.

Let $G(z)=\sum_{n=1}^\infty g_nz^n$ be the generating function for the sequence $\F{G}$. Note that $G(z)$ is holomorphic inside the unit disk. The recurrence relations in \eqref{recs} lead to a functional equation for $G(z)$.

\begin{lemma}\label{func} If $G(z)=\sum_n g_nz^n$, then $$G(z^2)=G(z)-\frac{z}{1+z^2}.$$
\end{lemma}

\begin{proof} This is directly given by the recurrences for $\F{G}$. We calculate that 
\begin{align*} G(z) &= \sum_n g_n z^n 
= \sum_n g_{2n} z^{2n}+\sum_n g_{4n+1} z^{4n+1}+\sum_n g_{4n+3} z^{4n+3}\\
&= \sum_n g_{n} z^{2n}+z\sum_n g_{4n+1} z^{4n}+z^3\sum_n g_{4n+3} z^{4n}\\
&= \sum_n g_{n} z^{2n}+z\sum_n g_{4n+1} z^{4n}-z^3\sum_n g_{4n+1} z^{4n}\\
&=G(z^2)+(z-z^3)\sum_n z^{4n}\\
&=G(z^2)+\frac{z-z^3}{1-z^4}.\end{align*} A little arithmetic and rearrangement gives the desired result.
\end{proof}

Using the functional equation from Lemma \ref{func}, we have $$G(z^{2^m})=G(z)-\sum_{k=0}^{m-1}\frac{z^{2^k}}{1+z^{2^{k+1}}}.$$ Denote the sum by $W_m(z)$; that is, $$W_m(z)=\sum_{k=0}^{m-1}\frac{z^{2^k}}{1+z^{2^{k+1}}},$$ so that \begin{equation}\label{GW} G(z^{2^m})=G(z)-W_m(z).\end{equation} 

\begin{proposition} If $|z|<1$, then the generating function of $\F{G}$ has the closed form  $$G(z)=\sum_{k=0}^{\infty}\frac{z^{2^k}}{1+z^{2^{k+1}}}.$$
\end{proposition}

\begin{proof} Take the limit as $m\to\infty$ in relation \eqref{GW}.
\end{proof}

Note that $$\gamma=G\left(\frac{1}{2}\right)=\lim_{m\to\infty}\left[G(2^{-2^m})+W_m\left(\frac{1}{2}\right)\right]=\sum_{k=0}^\infty\frac{1}{2^{2^k}+2^{-2^{k}}}.$$ 

\subsection{Transcendence of $G(z)$}

\begin{theorem}\label{transfunc} The function $G(z)$ is transcendental.
\end{theorem}

\begin{proof} Towards a contradiction, suppose that $G(z)$ is algebraic in $\B{C}[z]$; that is, there is an $n\geq 1$ and rational functions $q_0(z),q_1(z),\ldots,q_{n-1}(z)$ such that $$G(z)^n+q_{n-1}(z)G(z)^{n-1}+\cdots+q_0(z)=0\quad (|z|<1).$$

Choose $n$ minimally. By the functional equation of Lemma \ref{func} we obtain $$\sum_{k=0}^n q_k(z^2)\left[G(z)-\frac{z}{1+z^2}\right]^k=0\quad (|z|<1),$$ so that for $|z|<1$ \begin{align*}P(z):&=\sum_{k=0}^n q_k(z^2)[(1+z^2)G(z)-z]^k[1+z^2]^{n-k}\\
&=\sum_{k=0}^n q_k(z^2)[1+z^2]^{n-k}\sum_{j=0}^k{{k}\choose{j}}(1+z^2)^jG(z)^j(-z)^{k-j}=0.\end{align*} Thus \begin{equation}\label{Q} Q(z):=q_n(z)P(z)-(1+z^2)^nq_n(z^2)\sum_{k=0}^n q_k(z)G(z)^k=0.\end{equation} 

Inspection of $Q(z)$ gives the $k=n$ term as $$q_n(z)q_n(z^2)\left(\sum_{j=0}^n{{n}\choose{j}}(1+z^2)^jG(z)^j(-z)^{n-j}-(1+z^2)^nG(z)^n\right).$$ The coefficient of $G^n(x)$ is given when we set $j=n$ in the proceeding line, and is $$q_n(z)q_n(z^2)\left({{n}\choose{n}}(1+z^2)^nG(z)^n(-z)^{n-n}-(1+z^2)^nG(z)^n\right)=0.$$ Hence \eqref{Q} defines polynomials $h_1(z),\ldots,h_{n-1}(z)$ such that $$Q(z)=\sum_{k=0}^{n-1}h_k(z)G(z)^k=0.$$ The minimality of $n$ implies that $h_k(z)=0$ for $k=0,\ldots,n-1$.

Let us now determine $h_k(z)$ using the definition of $Q(z)$ from \eqref{Q}. We have \begin{align*}Q(z)&=\sum_{k=0}^{n-1}h_k(z)G(z)^k\\
&=\sum_{k=0}^n \Bigg\{\sum_{j=0}^k{{k}\choose{j}}q_n(z)q_k(z^2)(1+z^2)^{n-k+j}G(z)^j(-z)^{k-j}\\
&\qquad\qquad-(1+z^2)^nq_n(z^2)q_k(z)G(z)^k.\Bigg\}
\end{align*} From here we can read off the coefficient of $G(z)^m$ as
\begin{equation*}h_m(z)=\sum_{k=m}^n {{k}\choose{m}}q_k(z^2)(1+z^2)^{n-k+m}(-z)^{k-m}-(1+z^2)^nq_m(z),\end{equation*} recalling that $q_n(z)=1$.
Since $h_{n-1}(z)=0$, $$\sum_{k=n-1}^n {{k}\choose{n-1}}q_k(z^2)(1+z^2)^{n-k+n-1}(-z)^{k-(n-1)}=(1+z^2)^nq_{n-1}(z),$$ which leads to 
$$q_{n-1}(z^2)(1+z^2)^{n}-nz(1+z^2)^{n-1}=(1+z^2)^nq_{n-1}(z),$$ so that we may focus on the equation $$(1+z^2)q_{n-1}(z^2)-nz=(1+z^2)q_{n-1}(z).$$ Write $q_{n-1}(z)=\frac{a(z)}{b(z)}$, where $a(z)$ and $b(z)$ are polynomials with no common factors. Substituting and clearing denominators, we gain \begin{equation}\label{gcd}(1+z^2)a(z^2)b(z)-nzb(z)b(z^2)=(1+z^2)a(z)b(z^2).\end{equation} Using simple divisibility rules, \eqref{gcd} gives the following two conditions: $$ (i)\ \left.\frac{b(z)}{G(x)}\right| (1+z^2) \frac{b(z^2)}{G(x)}\qquad\mbox{and}\qquad (ii)\ \left.\frac{b(z^2)}{G(x)}\right| (1+z^2) \frac{b(z)}{G(x)},$$ where $G(x)=\gcd(b(z),b(z^2)).$ Recall that $(1+z^2)=(1+iz)(1-iz).$ 

A side note on determining properties of $b(z)$: condition $(ii)$ indicates the degree relationship, $2\deg b(z)\leq 2+\deg b(z)$, and equation \eqref{gcd} implies that $(z^2+1)|b(z)b(z^2)$ which gives $2\leq 3 \deg b(z).$ Together this yields a degree condition on $b(z)$ of $1\leq \deg b(z)\leq 2,$ since the degree must be a positive integer. In light of conditions $(i)$ and $(ii)$, we have $\deg b(z)=2.$

Now conditions $(i)$ and $(ii)$ imply that either \begin{equation}\label{1} (1+iz)\left|b(z) \right.\qquad\mbox{and}\qquad (1-iz)\left|b(z^2),\right. \end{equation} or \begin{equation}\label{2} (1-iz)\left|b(z) \right.\qquad\mbox{and}\qquad (1+iz)\left|b(z^2).\right.\end{equation} Given the above conditions, we have two options for $b(z)$: $$\mbox{condition \eqref{1}}\quad\Longrightarrow\quad b(z)=(z+1)(z-i),$$ or $$\mbox{condition \eqref{2}}\quad\Longrightarrow\quad b(z)=(z+1)(z+i).$$ 

Let us assume that condition \eqref{1} holds, that $b(z)=(z+1)(z-i)$. Then \eqref{gcd} becomes $$(1+z^2)(z+1)(z-i)a(z^2)-nz(z+1)(z-i)(z^2+1)(z^2-i)= (1+z^2)(z^2+1)(z^2-i)a(z).$$ Removing common factors, the last equation becomes \begin{equation}\label{a}(z+1)a(z^2)-nz(z+1)(z^2-i)= (z+i)(z^2-i)a(z).\end{equation} Equation \eqref{a} implies that $(z+1)|a(z)$, and so $(z^2+1)|a(z^2)$. Thus there exist $k(z),l(z)$ such that $$a(z)=k(z)(z+1)\quad\mbox{and}\quad a(z^2)=l(z)(z^2+1)=l(z)(z+i)(z-i).$$ Equation \eqref{a} becomes $$(z+1)l(z)(z+i)(z-i)-nz(z+1)(z^2-i)= (z+i)(z^2-i)k(z)(z+1),$$ implying that $$(z+i)|nz(z+1)(z^2-1),$$ which is not possible. Hence it must be the case that $b(z)\neq (z+1)(z-i)$, and so condition \eqref{1} is not possible.

It must now be the case that condition \eqref{2} holds, that $b(z)=(z+1)(z+i)$. Then \eqref{gcd} becomes $$(1+z^2)(z+1)(z+i)a(z^2)-nz(z+1)(z+i)(z^2+1)(z^2+i)= (1+z^2)(z^2+1)(z^2+i)a(z).$$ Removal of common factors yields \begin{equation}\label{a2}(z+1)a(z^2)-nz(z+1)(z^2+i)= (z-i)(z^2+i)a(z).\end{equation} Similar to the previous case, \eqref{a2} implies that $(z+1)|a(z)$, and so $(z^2+1)|a(z^2)$. Thus there exist $k(z),l(z)$ such that $$a(z)=k(z)(z+1)\quad\mbox{and}\quad a(z^2)=l(z)(z^2+1)=l(z)(z+i)(z-i).$$ Equation \eqref{a2} becomes $$(z+1)l(z)(z+i)(z-i)-nz(z+1)(z^2+i)= (z-i)(z^2+i)k(z)(z+1),$$ implying that $$(z+i)|nz(z+1)(z^2+i),$$ which is impossible. Since one of conditions \eqref{1} or \eqref{2} must hold, we arrive at our final contradiction, and the theorem is proved.
\end{proof}

\subsection{Transcendence of the Gaussian Liouville number $\gamma$}

We proceed to show that $\gamma$ is transcendental. We use a theorem of Mahler \cite{Mahl1}, as taken from Nishioka's book \cite{Nish1}. For the sake of completeness the proof of Mahler's theorem, taken again from \cite{Nish1}, is given (both the statement of Theorem \ref{Mahler} and its proof are taken verbatim from Nishioka's book). Here $\mathbf{I}$ is the set of algebraic integers over $\B{Q}$, $K$ is an algebraic number field, $\mathbf{I}_K=K\cap\mathbf{I},$ and $f(z)\in K[[z]]$ with radius of convergence $R>0$ satisfying the functional equation for an integer $d>1$, $$f(z^d)=\frac{\sum_{i=0}^m a_i(z)f(z)^i}{\sum_{i=0}^m b_i(z)f(z)^i},\qquad m<d,\ a_i(z),b_i(z)\in \mathbf{I}_K[z],$$ and $\gD(z):=\mbox{Res}(A,B)$ is the resultant of $A(u)=\sum_{i=0}^m a_i(z)u^i$ and $B(u)=\sum_{i=0}^m b_i(z)u^i$ as polynomials in $u$. Also, $$\overline{|\ga|}:=\max\{|\ga^\gs|:\gs\in\mbox{Aut}(\overline{\B{Q}}/\B{Q})\}\quad\mbox{and}\quad\mbox{den}(\ga):=\min\{d\in\B{Z}:d>0,d\ga\in\mathbf{I}\}.$$

\begin{theorem}[Mahler, \cite{Mahl1}]\label{Mahler} Assume that $f(z)$ is not algebraic over $K(z)$. If $\ga$ is an algebraic number with $0<|\ga|<\min\{1,R\}$ and $\gD(\ga^{d^k})\neq 0$ $(k\geq 0)$, then $f(\ga)$ is transcendental.
\end{theorem}

\begin{proof} Suppose that $f(\ga)$ is algebraic. We may assume $\ga,f(\ga)\in K$. Let $n$ be a positive integer. Then there are $n+1$ polynomials $P_0,P_1,\ldots,P_n\in \mathbf{I}_K[z]$ with degrees at most $n$ such that the auxiliary function $$E_n(z)=\sum_{j=0}^n P_j(z)f(z)^j=\sum_{h=0}^\infty b_h z^h$$ is not identically zero and all the coefficients $b_h$, with $h<n^2$, vanish. Since $f(z)$ is not algebraic over $K(z)$, $E_n(z)$ is not identically zero. Let $H$ be the least integer such that $b_H\neq 0$. Then $H>n^2$. Since $$\lim_{z\to 0}E_n(z)z^{-H}=b_H,$$ we have for any $k\geq c_1(n)$, \begin{equation}\label{1.2.2} 0\neq |E_n(\ga^{d^k})|\leq c_2(n)|\ga|^{d^kH}\leq c_2|\ga|^{d^kn^2}.\end{equation} There are polynomials $S(z,u),T(z,u)\in\mathbf{I}_K[z,u]$ such that $$\gD(z)=S(z,u)\sum_{i=0}^ma_i(z)u^i+T(z,u)\sum_{i=0}^m b_i(z)u^i.$$ Hence $$\gD(\ga)=S(\ga,f(\ga))\sum_{i=0}^ma_i(\ga)f(\ga)^i+T(\ga,f(\ga))\sum_{i=0}^m b_i(\ga)f(\ga)^i.$$ Suppose that $\sum_{i=0}^m b_i(\ga)f(\ga)^i=0$. Since $$\left(\sum_{i=0}^m b_i(\ga)f(\ga)^i\right)f(\ga^d)=\sum_{i=0}^m a_i(\ga)f(\ga)^i,$$ we get $\sum_{i=0}^m a_i(\ga)f(\ga)^i=0$ and so $\gD(\ga)=0$. This contradicts the assumption. Therefore $\sum_{i=0}^m b_i(\ga)f(\ga)^i\neq 0$ and $f(\ga^d)\in K$. Proceeding in this way, we see that $f(\ga^{d^k})\in K$ and therefore $E_n(\ga^{d^k})\in K$ $(k\geq 0)$. Define $Y_k$ $(k\geq 0)$ inductively as follows, \begin{align*} Y_1 &= \sum_{i=0}^m b_i(\ga)f(\ga)^i,\\
Y_{k+1} &= Y_k^m\sum_{i=0}^m b_i(\ga^{d^k})f(\ga^{d^k})^i,\quad k\geq 1.\end{align*} Then $Y_k\in K$ and $Y_k\neq 0$ $(k\geq 0)$. We estimate $\overline{|Y_k^n E_n(\ga^{d^k})|}$ and den$(Y_k^n E_n(\ga^{d^k})).$ Let $\deg_z(b)i(z)\leq l,\overline{|\ga|},\overline{|f(\ga)|}\leq c_3$ ($c_3>1$) and $D$ a positive integer such that $D\ga,Df(\ga)\in\mathbf{I}$. Then we have \begin{align*} \overline{|Y_1|} &= \overline{\left|\sum_{i=0}^m b_i(\ga)f(\ga)^i\right|}= \sum_{i=0}^m \overline{|b_i(\ga)|}\ \overline{|f(\ga)|}^i\leq c_4c_3^lc_3^m,\\
\overline{|Y_1f(\ga^d)|} &= \overline{\left|\sum_{i=0}^m a_i(\ga)f(\ga)^i\right|}= \sum_{i=0}^m \overline{|a_i(\ga)|}\ \overline{|f(\ga)|}^i\leq c_4c_3^lc_3^m\end{align*} and $$D^{l+m}Y_1,D^{l+m}Y_1f(\ga^d)\in\mathbf{I}.$$ Since $Y_2=Y_1^m\sum_{i=0}^mb_i(\ga^d)f(\ga^d)^i$ and $Y_2f(\ga^{d^k})=Y_1^m\sum_{i=0}^ma_i(\ga^d)f(\ga^d)^i$, we have $$\overline{|Y_2|},\ \overline{|Y_2f(\ga^{d^2})|}\leq (c_4c_3^{dl})(c_4c_3^{l+m})^m$$ and $$D^{dl}(D^{l+m})^mY_2,\ D^{dl}(D^{l+m})^mY_2f(\ga^{d^2})\in\mathbf{I}.$$ Proceeding in this way, we obtain $$\overline{|Y_k|},\ \overline{|Y_kf(\ga^{d^k})|}\leq c_4^{1+m+\cdots+m^{k-1}}(c_3^l)^{d^{k-1}+d^{k-2}m+\cdots+m^{k-1}}c_3^{m^k}$$ and \begin{align*} (D^l)^{d^{k-1}+d^{k-2}m+\cdots+m^{k-1}}D^{m^k}Y_k &\in\mathbf{I},\\
(D^l)^{d^{k-1}+d^{k-2}m+\cdots+m^{k-1}}D^{m^k}Y_kf(\ga^{d^k}) &\in\mathbf{I}. \end{align*} By the assumption $m<d$, we have $$d^{k-1}+d^{k-2}m+\cdots+m^{k-1}=d^{k-1}\left(1+\frac{m}{d}+\cdots+\left(\frac{m}{d}\right)^{k-1}\right)\leq c_5d^{k-1},$$ where we take a positive integer as $c_5$. Hence $$\overline{|Y_k|},\ \overline{|Y_kf(\ga^{d^k})|}\leq c_4^{c_5d^{k-1}}(c_3^l)^{c_5d^{k-1}}c_3^{d^k}\leq c_6^{d^k}$$ and $$D_0^{d^k}Y_k,\ D_0^{d^k}Y_kf(\ga^{d^k})\in\mathbf{I},\qquad D_0=D^{lc_5+1}.$$ Since $$Y_k^nE_n(\ga^{d^k})=\sum_{j=0}^n P_j(\ga^{d^k})Y_k^{n-j}\left(Y_kf(\ga^{d^k})\right)^j,$$ we obtain \begin{equation}\label{1.2.3}\overline{|Y_k^n E_n(\ga^{d^k})|}\leq c_7(n)c_3^{d^kn}c_6^{d^kn},\quad D_0^{2d^k}Y_k^nE_n(\ga^{d^k})\in\mathbf{I}.\end{equation} By \eqref{1.2.2}, \eqref{1.2.3} and the fundamental inequality, \begin{multline*} d^kn\log c_6+\log c_2(n)+d^kn^2\log |\ga|\geq \log|Y_k^nE_n(\ga^{d^k})|\\ \geq -2[K:\B{Q}] (\log c_7(n)+d^kn\log c_3c_6+2d^kn\log D_0),\end{multline*} for $k>c_1(n).$ Dividing both sides above by $d^k$ and letting $k$ tend to infinity, we have $$n\log c_6+n^2\log|\ga|\geq -2[K:\B{Q}](n\log c_3c_6+2n\log D_0).$$ Dividing both sides above by $n^2$ and letting $n$ tend to infinity, we have $\log|\ga|\geq 0$, a contradiction.
\end{proof}

Lemma \ref{func}, Theorem \ref{transfunc}, and Theorem \ref{Mahler} give our next theorem.

\begin{theorem} The Gaussian Liouville number $$\gamma=G\left(\frac{1}{2}\right)=\sum_{k=0}^\infty\frac{1}{2^{2^k}+2^{-2^{k}}}$$ is transcendental.
\end{theorem}

\begin{proof} Lemma \ref{func} gives the functional equation $$G(z^2)=\frac{(1+z^2)G(z)-z}{1+z^2},$$ so that, in the language of Theorem \ref{Mahler}, we have $$A(u)=(1+z^2)u-z\qquad\mbox{and}\qquad B(u)=1+z^2,$$ $m=1<2=d$, and $a_i(z),b_i(z)\in \mathbf{I}_K[z]$. Since $B(u)$ is a constant polynomial in $u$, $$\gD(z):=\mbox{Res}(A,B)=1+z^2.$$ Set $\ga=\frac{1}{2}$; it is immediate that $\ga=\frac{1}{2}$ is algebraic, $0<|\ga|=\frac{1}{2}<\min\{1,R\}$ ($R=1$), and $$\gD(\ga^{2^k})=\gD(2^{-2^k})=1+2^{-2^{k+1}}\neq 0\quad (k\geq 0).$$ Since $G(z)$ is not algebraic over $\B{C}[z]$ (as supplied by Theorem \ref{transfunc}), applying Theorem \ref{Mahler}, we have that $\gamma$ is transcendental.
\end{proof}

\section{Transcendence related to character--like functions}\label{gensection}

A {\em character--like function} $f$ associated to $p$ is a completely multiplicative function from $\B{N}$ to $\{-1,1\}$ defined by $f_1=1$, $f_{p}=\pm 1$ (your choice), and $f_{kp+i}=f_i$. As an example, the completely multiplicative function defined by $$f_n=\begin{cases} \pm 1 &\mbox{if } n=p\\ \leg{n}{p} &\mbox{if } p\nmid n,\end{cases}$$ where $\leg{n}{p}$ is the Legendre symbol modulo $p$, is a character--like function.  

If we let $f$ be a character--like function associated to $p$, then for $F(z)=\sum_n f_nz^n$ the generating function of the sequence $\F{F}:=(f_n)$ we have a lemma similar to the previous sections.

\begin{lemma}\label{genfunc} The generating function $F(z)$ of the sequence $\F{F}$ satisfies the functional equation $$F(z)=f_pF(z^p)+\frac{\Phi(z)}{1-z^p},$$ where $\Phi(z)=\sum_{i=1}^{p-1}f_i z^i.$
\end{lemma}

\begin{proof} We have \begin{align*} F(z) &= \sum_{k=0}^\infty \sum_{i=1}^{p-1} f_{pk+i}z^{pk+i}+\sum_{k=1}^\infty f_{pk}z^{pk}\\
&= \sum_{i=1}^{p-1} f_{i}z^i\sum_{k=0}^\infty z^{pk}+f_p\sum_{k=1}^\infty f_{k}z^{pk}\\
&= \frac{\sum_{i=1}^{p-1} f_{i}z^i}{1-z^p}+f_pF(z^p).\qedhere\end{align*}
\end{proof}

\begin{theorem}\label{gentransfunc} The function $F(z)$ is transcendental.
\end{theorem}

\begin{proof} We proceed as in the proof of Theorem \ref{transfunc}. Towards a contradiction, suppose that $F(z)$ is algebraic in $\B{C}[z]$; that is, there is an $n\geq 1$ and rational functions $q_0(z),q_1(z),\ldots,q_{n-1}(z)$ such that $$F(z)^n+q_{n-1}(z)F(z)^{n-1}+\cdots+q_0(z)=0\quad (|z|<1).$$

Choose $n$ minimally. The functional equation for $F(z)$ gives $$\sum_{k=0}^n q_k(z^2)\left[f_pF(z)+f_p\frac{\gF(z)}{z^p-1}\right]^k=0\quad (|z|<1),$$ where as before $\gF(z)=\sum_{i=1}^{p-1}f_i z^i.$ For $|z|<1$, we have \begin{align*}P(z):&=\sum_{k=0}^n q_k(z^p)f_p^k[(z^p-1)F(z)+\gF(z)]^k[1-z^p]^{n-k}\\
&=\sum_{k=0}^n q_k(z^p)f_p^k[z^p-1]^{n-k}\sum_{j=0}^k{{k}\choose{j}}(z^p-1)^jF(z)^j(\Phi(z))^{k-j}=0.\end{align*} Thus \begin{equation}\label{QG} Q(z):=q_n(z)P(z)-f_p^n(z^p-1)^nq_n(z^p)\sum_{k=0}^n q_k(z)F(z)^k=0.\end{equation} 

Inspection of $Q(z)$ gives the $k=n$ term as $$q_n(z)q_n(z^p)f_p^n \left(\sum_{j=0}^n{{n}\choose{j}}(z^p-1)^jF(z)^j(\Phi(z))^{n-j}-(z^p-1)^nF(z)^n\right).$$ The coefficient of $F^n(x)$ is given setting $j=n$ in the proceeding line, and is $$q_n(z)q_n(z^p)f_p^n\left({{n}\choose{n}}(z^p-1)^nF(z)^n(\Phi(z))^{n-n}-(z^p-1)^nF(z)^n\right)=0.$$ Hence \eqref{QG} defines polynomials $h_1(z),\ldots,h_{n-1}(z)$ such that $$Q(z)=\sum_{k=0}^{n-1}h_k(z)F(z)^k=0.$$ The minimality of $n$ gives that $h_k(z)=0$ for $k=0,\ldots,n-1$.

Let us now determine $h_k(z)$ using the definition of $Q(z)$ from \eqref{QG}. We have \begin{align*}Q(z)&=\sum_{k=0}^{n-1}h_k(z)F(z)^k\\
&=\sum_{k=0}^n \Bigg\{\sum_{j=0}^k{{k}\choose{j}}q_n(z)q_k(z^p)f_p^k(z^p-1)^{n-k+j}F(z)^j(\Phi(z))^{k-j}\\
&\qquad\qquad-f_p^n(z^p-1)^nq_n(z^p)q_k(z)F(z)^k.\Bigg\}
\end{align*} From here we can read off the coefficient of $F(z)^m$ as
\begin{equation*}h_m(z)=\sum_{k=m}^n {{k}\choose{m}}q_k(z^p)f_p^k(z^p-1)^{n-k+m}(\Phi(z))^{k-m}-f_p^n(z^p-1)^nq_m(z),\end{equation*} recalling that $q_n(z)=1$.

Since $h_{n-1}(z)=0$, we have $$\sum_{k=n-1}^n {{k}\choose{n-1}}q_k(z^p)f_p^k(z^p-1)^{n-k+n-1}(\Phi(z))^{k-(n-1)}= f_p^n(z^p-1)^nq_{n-1}(z),$$ which yields 
$$q_{n-1}(z^p)f_p^{n-1}(z^p-1)^{n}+n f_p^n \gF(z)(z^p-1)^{n-1}= f_p^n(z^p-1)^nq_{n-1}(z),$$ so that we may focus on the equation $$f_p(z^p-1)q_{n-1}(z^p)+n \gF(z)= (z^p-1)q_{n-1}(z).$$ 

Write $q_{n-1}(z)=\frac{a(z)}{b(z)}$, where $a(z)$ and $b(z)$ are polynomials with no common factors. Substituting and clearing denominators, we gain \begin{equation}\label{gengcd}f_p(z^p-1)a(z^p)b(z)+n \gF(z)b(z)b(z^p)= (z^p-1)a(z)b(z^p).\end{equation} Again, using divisibility rules, \eqref{gcd} gives the following two conditions: $$ (i)\ \left.b(z)\right| (z^p-1) b(z^p)\qquad\mbox{and}\qquad (ii)\ \left.b(z^p)\right| (z^p-1) b(z).$$ 

Condition $(ii)$ indicates the degree relationship, that $p\deg b(z)\leq p+\deg b(z)$, which gives (since $p>2$) $$\deg b(z)\leq \frac{p}{p-1}<2.$$ We have two choices, either $\deg b(z)=0$ or $\deg b(z)=1$.

If $\deg b(z)=0$, then $b(z)\in\B{C}$, say $b(z)=b$ (which also gives $b(z^2)=b$). Equation \eqref{gengcd} becomes \begin{equation*}f_p(z^p-1)a(z^p)b+n \gF(z)b^2= (z^p-1)a(z)b.\end{equation*} This equation gives $(z^p-1)|\gF(z)$ which implies that $\deg \gF(z)\geq p,$ an  impossibility since the definition of $\gF(z)$ gives $\deg \gF(z)=p-1$. Thus it must be the case that $\deg b(z)=1.$

If $\deg b(z)=1$, then $b(z)=z-\gb$ (the constant can be given to $a(z)$). In this case, \eqref{gengcd} becomes \begin{equation}\label{gengcd1} f_p(z^p-1)a(z^p)(z-\gb)+n \gF(z)(z-\gb)(z^p-\gb)= (z^p-1)a(z)(z^p-\gb).\end{equation} Since $\gcd(a(z),b(z))=1$, equation \eqref{gengcd1} gives that $$(z^p-\gb)|(z^p-1)(z-\gb).$$ This implies that $\gb=1$. Thus \eqref{gengcd1} becomes \begin{equation*} f_p(z^p-1)a(z^p)(z-1)+n \gF(z)(z-1)(z^p-1)= a(z)(z^p-1)^2.\end{equation*} With repeated factors deleted this is \begin{equation*} f_pa(z^p)(z-1)-n \gF(z)(z-1)= a(z)(z^p-1);\end{equation*} that is, \begin{equation}\label{gengcd12} (z-1)(f_pa(z^p)+n \gF(z))= a(z)(z^p-1).\end{equation} Comparing the degrees of each side gives $$p \deg a(z)+1=\deg a(z)+p,$$ so that $\deg a(z)=1$. Write $a(z)=\gd(z-\ga).$ Substituting this into \eqref{gengcd12} yields \begin{equation*} (z-1)(f_p\gd(z^p-\ga)+n \gF(z))= \gd(z-\ga)(z^p-1).\end{equation*} Since $\deg \gF(z)=p-1$ comparing lead coefficients gives $f_p\gd=\gd$, so that $f_p=1$. (In the case that $f_p$ is defined as $-1$ this gives a contradiction.)

Removing the $(z-1)$ factor gives \begin{equation}\label{ga}\gd(z^p-\ga)+n \gF(z)= \gd(z-\ga)\sum_{m=0}^{p-1}z^m.\end{equation} Note that the right--hand side and left--hand side of \eqref{ga} are equal as polynomials; that is, term by term. If we expand out \eqref{ga} we have \begin{equation*}\label{coef} \gd z^p+n\sum_{m=1}^{p-1}f_m z^m -\gd\ga= \gd z^p-\gd(1+\ga)\sum_{m=1}^{p-1}z^m-\gd\ga.\end{equation*} Hence for every $m=1,2,\ldots,p-1,$ we have that \begin{equation}\label{nda}nf_m=-\gd(1+\ga),\end{equation} and specifically when $m=1$ this tells us that $n=-\gd(1+\ga)$ since $f_1=1$. Thus \eqref{nda}, with the substitution $n=-\gd(1+\ga)$, gives $f_m=1$ for each $m=1,2,\ldots,p-1$. This is impossible for an odd prime $p$, and the theorem is proven.
\end{proof}

Using the functional equation, we yield $$F(z^p)=f_pF(z)-f_p\frac{\Phi(z)}{1-z^p}.$$ We can build a series for the number as before. The functional equation gives \begin{align*}F(z^{p^{m}})&=f_pF(z^{p^{m-1}})-f_p\frac{\Phi(z^{p^{m-1}})}{1-z^{p^{m}}}\\
&=f_p^2F(z^{p^{m-2}})-f_p^2\frac{\Phi(z^{p^{m-2}})}{1-z^{p^{m-1}}}-f_p\frac{\Phi(z^{p^{m-1}})}{1-z^{p^{m}}}\\
&=f_p^mF(z)-\sum_{k=1}^mf_p^k\frac{\Phi(z^{p^{m-k}})}{1-z^{p^{m-k+1}}},
\end{align*} which when rearranged leads to $$F(z)=f_p^mF(z^{p^{m}})+\sum_{k=1}^m f_p^{m-k}\frac{\Phi(z^{p^{m-k}})}{1-z^{p^{m-k+1}}}.$$ We change the index and set $$V_m(z):=\sum_{k=0}^{m-1} f_p^{k}\frac{\Phi(z^{p^{k}})}{1-z^{p^{k+1}}}$$
to give $$F(z)=f_p^mF(z^{p^{m}})+V_m(z).$$

\begin{proposition}\label{Fseries} If $|z|<1$, then the generating function has the closed form  $$F(z)=\sum_{k=0}^{\infty} f_p^{k}\frac{\Phi(z^{p^{k}})}{1-z^{p^{k+1}}}.$$
\end{proposition}

\begin{proof} Take the limit as $m\to\infty$ in the equation $F(z)=f_p^mF(z^{p^{m}})+V_m(z)$.
\end{proof}

At $z=1/2$ we have $$F\left(\frac{1}{2}\right)=\lim_{m\to\infty}\left[f_p^mF(2^{-p^{m}})+V_m\left(\frac{1}{2}\right)\right]=\sum_{k=0}^\infty \frac{f_p^k\Phi(2^{-p^{k}})}{1-2^{-p^{k+1}}}.$$

\begin{theorem} For each odd prime $p$, the number $\gf_p:=F\left(\frac{1}{2}\right)$ is transcendental.
\end{theorem}

\begin{proof} Lemma \ref{genfunc} gives \begin{equation*} F(z^p)=\frac{f_p(1-z^p)F(z)-f_p\gF(z)}{1-z^p},\end{equation*} where $\gF(z)=\sum_{i=1}^{p-1}f_i z^i$. Similar to the specific case of the previous section, using the language of Mahler's theorem, we have $$A(u)=f_p(1-z^p)u-f_p\gF(z)\qquad\mbox{and}\qquad B(u)=1-z^p,$$ $m=1<p=d$, and $a_i(z),b_i(z)\in \mathbf{I}_\B{C}[z]$. Again, $B(u)$ is a constant polynomial in $u$, so that $$\gD(z):=\mbox{Res}(A,B)=1-z^p.$$ Set $\ga=\frac{1}{2}$; $\ga=\frac{1}{2}$ is algebraic, $0<|\ga|=\frac{1}{2}<\min\{1,R\}$ ($R=1$), and $$\gD(\ga^{p^k})=\gD(2^{-p^k})=1-2^{-p^{k+1}}\neq 0\quad (k\geq 0).$$ Theorem \ref{gentransfunc} gives that $F(z)$ is not algebraic over $\B{C}[z]$ and we may apply Theorem \ref{Mahler} to give the desired result.
\end{proof}


\remark{Mahler's theorem tells us that the values of the functions $G(z)$ and $F(z)$ are transcendental for any algebraic value of $z$ within their radii of convergence. The special value of $z=\frac{1}{2}$ is focused on only for its relation to the sequences $\F{G}$ and $\F{F}$.}

\bibliographystyle{amsplain}

\end{document}